\documentstyle{amsart}
\input{bull-art}

\begin{document}
\def\currentvolume{30}
\def\currentissue{2}
\def\currentyear{1994}
\def\currentmonth{April}
\def\copyrightyear{1994}
\def\currentpages{161-177}
\title{On proof and progress in mathematics}
\author{William P. Thurston}
\date{October 26, 1993}
\subjclass{Primary 01A80}
\maketitle

This essay on the nature of proof and progress in 
mathematics was 
stimulated by the article of Jaffe and Quinn, 
``Theoretical Mathematics: 
Toward a cultural synthesis of mathematics and theoretical 
physics".
Their article raises interesting issues that 
mathematicians should pay more
attention to, but it also perpetuates some widely held
beliefs and attitudes that need to be questioned and 
examined.

The article had one paragraph portraying some of my work 
in a way
that diverges from my experience, and it also diverges 
from the observations
of people in the field whom I've discussed it with as a 
reality check.

After some reflection, it seemed to me that what Jaffe and 
Quinn
wrote was an example of the phenomenon that people see 
what they are
tuned to see. Their portrayal of my work resulted from
projecting the sociology of mathematics onto a 
one-dimensional
scale (speculation versus rigor) that ignores many basic 
phenomena.

Responses to the Jaffe-Quinn article have been invited 
from a number
of mathematicians, and I expect it to receive plenty of 
specific analysis
and criticism from others.  Therefore, I will concentrate 
in this essay
on the positive rather than on the contranegative.
I will describe my view of the process of mathematics, 
referring
only occasionally to Jaffe and Quinn by way of comparison.

\bigskip

In attempting to peel back layers of assumptions, it is 
important to
try to begin with the right questions:

\section{What is it that mathematicians accomplish?}

There are many issues buried in this question, which I 
have tried to phrase in
a way that does not presuppose the nature of the answer.

It would not be good to start, for example, with the 
question
\begin{quotation}
How do mathematicians prove theorems?
\end{quotation}  
This question introduces an interesting topic, but 
to start with it would be to project two hidden assumptions:
\begin{enumerate}  
\item that there is uniform, objective and firmly 
established theory and practice of mathematical proof, and
\item that progress made by mathematicians consists of 
proving 
theorems. \end{enumerate}
It is worthwhile to examine these hypotheses, rather than 
to accept
them as obvious and proceed from there.

The question is not even
\begin{quote}
How do mathematicians make progress in mathematics?
\end{quote} 
Rather, as a more explicit (and leading) form of the 
question, I prefer
\begin{quote}
How do mathematicians advance human understanding of 
mathematics?
\end{quote}

This question brings to the fore something that is
fundamental and pervasive: that what we are doing is finding
ways for {\it people} to understand and think about 
mathematics.

The rapid advance of computers has helped 
dramatize this point, because computers and people are 
very different.
For instance, when Appel and Haken completed a proof of 
the 4-color 
map theorem using a massive automatic computation, it 
evoked much controversy.
I interpret the controversy as having
little to do with doubt people had as to the veracity 
of the theorem or the correctness of the proof.  Rather, 
it reflected a 
continuing desire for {\it human understanding} of a 
proof, in addition
to knowledge that the theorem is true.

On a more everyday level, it is common for people first 
starting to 
grapple with computers to make large-scale computations
of things they  might have done on a smaller scale by hand.
They might print out a table of the first 10,000 primes,
only to find that their printout isn't 
something they really wanted after all.  They discover by 
this kind of
experience that what they really want is usually not some
collection of ``answers"---what they want is {\it 
understanding}. 

\bigskip

It may sound almost circular to say that what
mathematicians are accomplishing is
to advance human understanding of mathematics.
I will not try to resolve this by discussing what 
mathematics is,
because it would take us far afield. Mathematicians 
generally feel that
they know what mathematics is, but find it difficult to 
give a good direct
definition. It is interesting to try. For me, ``the theory 
of formal
patterns'' has come the closest, but to discuss this would 
be a whole essay in
itself.

Could the difficulty in giving a good direct definition of 
mathematics
be an essential one, indicating that mathematics has an 
essential
recursive quality? Along these lines we might say that 
mathematics
is the smallest subject satisfying the following:
\begin{itemize}
\item Mathematics includes the natural numbers and plane 
and solid geometry.
\item Mathematics is that which mathematicians study.
\item Mathematicians are those humans who advance human 
understanding of
mathematics.
\end{itemize}
In other words, as mathematics advances, we incorporate it 
into our thinking.
As our thinking becomes more sophisticated, we generate 
new mathematical
concepts and new mathematical structures: the subject 
matter of
mathematics changes to reflect how we think.

\medskip
If what we are doing is constructing better ways of 
thinking, then
psychological and social dimensions are essential to a 
good model
for mathematical progress.  These dimensions are absent 
from the popular
model.  In caricature, the popular model
holds that
\begin{itemize}
\item[\bf D.]
 mathematicians start from a few basic mathematical 
structures
and a collection of axioms ``given" about these 
structures, that 
\item[\bf T.]
there are various important questions 
to be answered about these structures that can be stated 
as formal mathematical
propositions, and
\item[\bf P.]
the task of the 
mathematician is to seek a deductive pathway from the 
axioms to 
the propositions or to their denials.
\end{itemize}
 We might call this the
definition-theorem-proof (DTP) model of mathematics.

A clear difficulty with the DTP model is that it doesn't 
explain the source of
the questions.  Jaffe and Quinn discuss
speculation (which they inappropriately label
``theoretical mathematics'') as an 
important additional ingredient.  Speculation
consists of making conjectures, 
raising questions, and making intelligent guesses and 
heuristic arguments
about what is probably true. 

Jaffe and Quinn's DSTP model still fails to address some 
basic issues.
We are not trying to meet some
abstract production quota of definitions, theorems and 
proofs.
The measure of our success is whether what we do enables 
{\it people} to
understand and think more clearly and effectively about 
mathematics.

Therefore, we need to ask ourselves:

\section{How do people understand mathematics?}

This is a very hard question. Understanding is 
an individual and internal matter that is hard to be fully 
aware of, hard to
understand and often hard to communicate.   We can only 
touch on it
lightly here.

People have very different ways of understanding 
particular pieces of mathematics. To illustrate this, it 
is best to
take an example that practicing mathematicians understand 
in multiple 
ways, but that we see our
students struggling with.  The derivative of a function
fits well.  The derivative can be thought of as:
\begin{enumerate}
\item Infinitesimal:  the ratio of the infinitesimal
change in the value of a function to the infinitesimal 
change in a 
function.
\item  Symbolic:
the derivative of $x^n$ is $n x^{n-1}$, the derivative of 
$\sin(x)$ is 
$\cos(x)$, the derivative of $f \circ g$ is $f'\circ g * 
g'$, {\it etc.}
\item Logical: $f'(x) = d$ if and only if for every 
$\epsilon$ there is
a $\delta$ such that when $0 < \left | \Delta x \right | < 
\delta$,
$$ \left | {{f(x + \Delta x) - f(x) } \over {\Delta x}} - 
d \right | < \delta. $$
\item  Geometric: the derivative is the slope of a line 
tangent to the 
graph of the function, if the graph has a tangent.
\item  Rate:  the instantaneous speed of $f(t)$, when $t$ 
is time.
\item   Approximation: The derivative of a function is the 
best linear
approximation to the function near a point.
\item  Microscopic:  The derivative of a function is the 
limit of what you
get by looking at it under a microscope of higher and 
higher power.
\end{enumerate}

This is a list of different
ways of {\it thinking about} or {\it conceiving of} the 
derivative,
rather than a list of different {\it logical definitions}. 
Unless great efforts are made to maintain the tone and 
flavor of
the original human insights, the differences start to 
evaporate as
soon as the mental concepts are translated into 
precise, formal and explicit definitions.

I can remember absorbing each of these concepts as 
something new and interesting, and spending a good deal of 
mental
time and effort digesting and practicing with each,
reconciling it with the others.  I also remember coming 
back to
revisit these different concepts later with added meaning 
and understanding.

The list continues; there is no reason for it ever to stop.
A sample entry further down the list may help illustrate 
this.  We may
think we know all there is to say about a certain subject,
but new insights are around the corner. Furthermore, one 
person's clear
mental image is another person's intimidation:
\begin{enumerate}
\item[37.]
The derivative of a real-valued function $f$ in a domain $D$
is the Lagrangian section of the cotangent bundle $T^*(D)$
that gives the connection form for the unique flat 
connection on the trivial
$\bold R$-bundle $D \times \bold R$ for which the graph of 
$f$ is parallel.
\end{enumerate}

These differences are not just a curiosity.  Human 
thinking and 
understanding do not work on a single track, like a 
computer with a 
single central processing unit.  Our brains and minds seem 
to be
organized into a variety of separate, powerful facilities.
These facilities work together loosely, ``talking" to each 
other at
high levels rather than at low levels of organization.

Here are some major divisions that are important for 
mathematical thinking: 
\begin{enumerate}
\item Human language.  We have powerful special-purpose
facilities for speaking and understanding human language, 
which also tie
in to reading and writing.  Our linguistic facility is an 
important
tool for thinking, not just for communication.  A crude
example is the quadratic formula which people may remember 
as a little chant,
``ex equals minus bee plus or minus the square root of bee 
squared
minus four ay see all over two ay."   The
mathematical language of symbols is closely tied to our 
human language facility.
The fragment of mathematical symbolese available to most 
calculus students
has only one verb, ``\<$=$''.  That's why students use it
when they're in need of a verb.  Almost anyone who has 
taught calculus
in the U.S. has seen students instinctively write ``\<$x^3 
= 3 x^2$'' and the
like.
\item  Vision, spatial sense, kinesthetic (motion) sense.
People have very powerful facilities for taking in 
information visually
or kinesthetically, and thinking with their spatial sense. 
 On the other
hand, they do not have a very good built-in facility for 
inverse vision,
that is, turning an internal spatial understanding back 
into a two-dimensional
image.  Consequently, mathematicians usually have fewer 
and poorer
figures in their papers and books than in their heads.

An interesting phenomenon in spatial thinking is that 
scale makes
a big difference.  We can think about little objects in 
our hands, or we can
think of bigger human-sized structures that we scan, or we 
can think of
spatial structures that encompass us and that we move 
around in. 
We tend to think more effectively with spatial imagery on 
a larger scale:
it's as if our brains take larger things more seriously 
and can
devote more resources to them.

\item  Logic and deduction.  We have some built-in ways of
reasoning and putting things together associated with how 
we make logical
deductions:  cause and effect (related to implication),
contradiction or negation, {\it etc.} 

Mathematicians apparently don't generally
rely on the formal rules of deduction as they are 
thinking.  Rather, they hold
a fair bit of logical structure of a proof in their heads, 
breaking 
proofs into intermediate results so that they don't have 
to hold too much 
logic at once.  In fact, it is common for excellent 
mathematicians not even to
know the standard
formal usage of quantifiers (for all and there exists), 
yet all
mathematicians certainly perform the reasoning that they 
encode.

It's interesting that although ``or'', ``and''
and ``implies'' have identical formal usage, we think of 
``or'' and ``and''
as conjunctions and ``implies'' as a verb.
\item  Intuition, association, metaphor.  People have 
amazing 
facilities for sensing something without knowing where it 
comes 
from (intuition); for sensing that some phenomenon or 
situation or 
object is like something else (association); and for 
building and testing connections 
and comparisons, holding two things in mind at the same 
time (metaphor).
These facilities are quite important for mathematics.  
Personally, I put a lot of effort into ``listening" to my 
intuitions and 
associations, and building them into metaphors and 
connections. This involves
a kind of simultaneous quieting and focusing of
my mind.  Words, logic, and detailed pictures rattling 
around can inhibit
intuitions and associations.
\item Stimulus-response.  This is often emphasized in 
schools; for instance, if you see $3927 \times 253$,  you 
write one number above the other and draw a line 
underneath, {\it etc.}
This is also important for research mathematics: seeing a 
diagram of a 
knot, I might write down a presentation for the 
fundamental group of its 
complement by a procedure that is similar in feel to the 
multiplication
algorithm.
\item  Process and time.  We have a facility for thinking 
about
processes or sequences of actions that can often be used 
to good effect 
in mathematical reasoning. One way to think of a function 
is as an 
action, a process, that takes the domain to the range.  
This is particularly
valuable when composing functions.  Another use of this 
facility is in remembering proofs:  people often remember 
a proof as a process
consisting of several steps. In topology, the notion of a 
homotopy is 
most often thought of as a process taking time.   
Mathematically, time is 
no different from one more spatial dimension, but since 
humans interact 
with it in a quite different way, it is psychologically 
very different.
\end{enumerate}

\section{How is mathematical understanding communicated?}

The transfer of understanding from one person
to another is not automatic. It is hard and tricky.
Therefore, to analyze human understanding of mathematics, 
it is
important to consider {\bf who} understands {\bf what}, 
and {\bf when}.

Mathematicians have 
developed habits of communication that are often 
dysfunctional. Organizers
of colloquium talks everywhere exhort speakers to explain 
things in 
elementary terms.  Nonetheless, most of the audience at an 
average 
colloquium talk gets little of value from it.  
Perhaps they are lost within the first
5 minutes, yet sit silently through the remaining 55 
minutes. Or perhaps they
quickly lose interest because the speaker plunges
into\ technical details
without presenting any reason to investigate them.   
At the end of the talk, the few mathematicians who are 
close to the field
of the speaker ask a question or two to avoid 
embarrassment.  

This pattern is similar to what often holds in 
classrooms, where we go through the motions of 
saying for the record what we think the students ``ought" 
to learn,
while the students are trying to grapple with the more 
fundamental issues of
learning our language and guessing at our mental models.   
Books compensate
by giving samples of how to solve every type of homework 
problem.
Professors compensate by giving homework and tests that 
are much easier than
the material ``covered" in the course, and then grading 
the homework
and tests on a scale that requires little understanding.
We assume that the problem is with the students rather 
than with
communication: that the students either just don't have 
what it takes,
or else just don't care.

Outsiders are amazed at this phenomenon, but within the
mathematical community, we dismiss it with shrugs.

\smallskip

Much of the difficulty has to do with 
the language and culture of mathematics, which is divided 
into 
subfields.   Basic concepts used every day 
within one subfield are often foreign to another subfield.
Mathematicians give up on trying to understand the basic 
concepts even 
from neighboring subfields, unless they were clued in as 
graduate 
students.

In contrast, communication works very well within the 
subfields
of mathematics. Within a subfield, people develop a body
of common knowledge and known techniques.   By informal 
contact, people 
learn to understand and copy each other's ways of thinking,
so that ideas can be explained clearly and easily.

Mathematical
knowledge can be transmitted amazingly fast within a 
subfield.  When a 
significant theorem is proved, it often (but not always) 
happens that
the solution can be communicated in a matter of minutes 
from one person
to another within the subfield.  The same proof would be 
communicated and
generally understood in an hour talk to members of the 
subfield.   It
would be the subject of a 15- or 20-page paper, which 
could be read and
understood in a few hours or perhaps days by members of 
the subfield.

\medskip
Why is there such a big expansion from the informal 
discussion to the
talk to the paper?  One-on-one, people use wide channels 
of communication
that go far beyond formal mathematical language.   They 
use gestures, they 
draw pictures and diagrams, they make sound effects and 
use body language.
Communication is more likely to be two-way, so that people 
can concentrate on
what needs the most attention.  With these channels of 
communication,
they are in a much better position 
to convey what's going on, not just in their logical and 
linguistic 
facilities, but in their other mental facilities as well.

In talks, people are more inhibited and more formal.   
Mathematical
audiences are often not very good at asking the questions 
that are
on most people's minds, and speakers often have an 
unrealistic preset outline
that inhibits them from addressing questions even when 
they are asked.

In papers, people 
are still more formal.  Writers translate their ideas into 
symbols and logic, and readers try to translate back.

\medskip

Why is there such a discrepancy between communication 
within a subfield 
and communication outside of subfields, not to mention 
communication
outside mathematics?

Mathematics in some sense has a common language: a 
language of symbols, technical definitions, computations, 
and logic. 
This language efficiently conveys some, but not all, modes 
of mathematical
thinking.  Mathematicians learn to translate certain 
things almost
unconsciously from one mental mode to the other, so that 
some statements
quickly become clear.
Different mathematicians study papers in different ways, 
but when
I read a mathematical paper in a field in which I'm 
conversant, I concentrate
on the thoughts that are between the lines.  I might look 
over several
paragraphs or strings of equations
and think to myself ``Oh yeah, they're
putting in enough rigamarole to carry such-and-such idea.'' 
When the idea is clear, the formal setup is usually 
unnecessary and
redundant---I often feel that I could write it out myself 
more easily
than figuring out what the authors actually wrote.   It's 
like
a new toaster that comes with a 16-page manual. If you 
already understand
toasters and if the toaster looks like previous toasters 
you've encountered,
you might just plug it in and see if it works, rather than 
first
reading all the details in the manual.

People familiar with ways of doing things in a subfield
recognize various patterns of statements or formulas as 
idioms or
circumlocution for certain concepts or mental images.  But 
to people not
already familiar with what's going on the same patterns 
are not
very illuminating; they are often even misleading.
The language is not alive except to those who use it.

\bigskip

I'd like to make an important remark here: there are some
mathematicians who are conversant with the ways of 
thinking in more than
one subfield, sometimes in quite a number of subfields.
Some mathematicians learn the jargon of several
subfields as graduate students, some people are just quick 
at
picking up foreign mathematical language and culture, and 
some
people are in mathematical centers where they are exposed 
to many
subfields.  People who are comfortable in more than one 
subfield can
often have a very positive influence, serving as bridges, 
and helping
different groups of mathematicians learn from each other.
But people knowledgeable in multiple fields can also have
a negative effect, by intimidating others,
and by helping to validate and maintain the whole system of
generally poor communication.  For example, one effect often
takes place during colloquium talks, where one or two widely
knowledgeable people sitting in the front row may serve as 
the speaker's
mental guide to the audience.

\bigskip

There is another effect caused by the big differences 
between how we
think about mathematics and how we write it.   A group
of mathematicians interacting with each other can keep a 
collection of
mathematical ideas alive for a period of years, even 
though the recorded
version of their mathematical work differs from their actual
thinking, having much greater emphasis on language,
symbols, logic and formalism.
But as new batches of mathematicians learn about the
subject they tend to interpret what they read and
hear more literally, so that the more easily recorded and 
communicated
formalism and machinery tend to gradually take over from 
other modes of
thinking.

There are two counters to this trend, so that mathematics
does not become entirely mired down in formalism.  First, 
younger
generations of mathematicians are continually discovering 
and rediscovering
insights on their own, thus reinjecting diverse modes of 
human thought
into mathematics.

Second, mathematicians sometimes invent names
and hit on unifying definitions that replace technical
circumlocutions and give good handles for insights.  
Names like ``group" to replace ``a system of substitutions 
satisfying \dots", 
and ``manifold" to replace 
\begin{quote}
We can't give coordinates to
parametrize all the solutions to our equations 
simultaneously,
but in the neighborhood of any particular solution we can 
introduce
coordinates
\begin{equation*}
\begin{aligned}
&(f_1(u_1,u_2,u_3), f_2(u_1,u_2,u_3), f_3(u_1,u_2,u_3), 
f_4(u_1,u_2,u_3),\\
&f_5(u_1,u_2,u_3))\end{aligned}\end{equation*}
where at least one of the ten determinants 
$$ \dots\! \text{[ten \!$3\!\times\!3$ \!determinants 
of matrices of partial derivatives]} \!\dots 
$$
is not zero
\end{quote}
may or may not have represented advances in insight among 
experts,
but they greatly facilitate the {\it communication} of 
insights.

\bigskip

We mathematicians need to put far greater effort into 
communicating mathematical {\it ideas}.
To accomplish this, we need to pay much more attention to 
communicating
not just our definitions, theorems, and proofs, but also
our ways of thinking.  We need to appreciate the value of 
different ways
of thinking about the same mathematical structure.

We need to focus far more
energy on understanding and explaining
the basic mental infrastructure of mathematics---with 
consequently
less energy on the most recent results.   This entails 
developing
mathematical language that is effective for the radical 
purpose of
conveying ideas to people who don't already know them.

Part of this communication is through proofs.

\section{What is a proof?}

When I started as a graduate student at Berkeley, I had 
trouble imagining
how I could ``prove'' a new and interesting mathematical 
theorem.
I didn't really understand what a ``proof'' was.

By going to seminars, reading papers, and talking to other 
graduate students,
I gradually began to catch on.  Within any field,
there are certain theorems and certain techniques that are 
generally
known and generally accepted.    When you write a paper, 
you refer
to these without proof.   You look at other papers in the 
field, and you
see what facts they quote without proof, and what they 
cite in their
bibliography.  You learn from other people some idea of 
the proofs.
Then you're free to quote the same theorem and
cite the same citations.  You don't necessarily have to 
read the full
papers or books that are in your bibliography.
Many of the things that are generally known are things for 
which
there may be no known written source.  As long as people 
in the field are
comfortable that the idea works, it doesn't need to have a 
formal written
source.

At first I was highly suspicious of this process.   I 
would doubt whether
a certain idea was really established.  But I found that I 
could ask
people, and they could produce explanations and proofs,
or else refer me to other people or to written sources 
that would give
explanations and proofs.  There were
published theorems that were generally known to be false, 
or where the
proofs were generally known to be incomplete.
Mathematical knowledge
and understanding were embedded in the minds and in the 
social fabric of the
community of people thinking about a particular topic.
This knowledge was supported by written documents, but the 
written documents
were not really primary.

I think this pattern varies quite a bit from field to 
field.   I was interested
in geometric areas of mathematics, where it is often 
pretty hard to have
a document that reflects well the way people actually 
think.   In more
algebraic or symbolic fields, this is not necessarily so, 
and I have the
impression that in some areas documents are much closer to 
carrying
the life of the field.   But in any field, there is a 
strong social standard
of validity and truth.   Andrew Wiles's proof of Fermat's 
Last Theorem is a
good illustration of this, in a field which is very 
algebraic. The experts
quickly came to believe that his proof was basically 
correct on the basis of
high-level ideas, long before details could be checked. 
This proof will receive
a great deal of scrutiny and checking compared to most 
mathematical
proofs; but no matter how the process of verification 
plays out, it helps
illustrate how mathematics evolves by rather organic 
psychological and
social processes.

When people are doing mathematics, the flow of
ideas and the social standard of validity is much more 
reliable than
formal documents.  People are usually not very good in 
checking
{\it formal correctness} of proofs, but they are quite 
good at
detecting potential weaknesses or flaws in proofs.

To avoid misinterpretation, I'd like to emphasize two 
things I am
{\it not} saying. First, I am {\it not} advocating any 
weakening
of our community standard of proof; I am trying to 
describe how the process
really works. Careful proofs that will stand up to 
scrutiny are very
important. I think the process of proof on the whole works 
pretty well
in the mathematical community. The kind of change I would 
advocate is that
mathematicians take more care with their proofs, making 
them really
clear and as simple as possible so that if any weakness is 
present
it will be easy to detect. Second, I am {\it not} 
criticizing the
mathematical study of formal proofs, nor am I criticizing 
people
who put energy into making mathematical arguments more 
explicit
and more formal. These are both useful activities that 
shed new
insights on mathematics.

\bigskip

I have spent a fair amount of effort during periods of my 
career
exploring mathematical questions by computer.  In view of 
that experience,
I was astonished to see the statement of Jaffe and Quinn 
that
mathematics is extremely slow and arduous, and that it is 
arguably
the most disciplined of all human activities.    The 
standard of
correctness and completeness necessary to get a computer 
program
to work at all is a couple of orders of magnitude higher 
than the
mathematical community's standard of valid proofs.
Nonetheless, large computer programs, even when they have 
been very carefully
written and very carefully tested, always seem to have bugs.

I think that mathematics is one of the most intellectually
gratifying of human activities. Because we have a high 
standard
for clear and convincing thinking and because we place a 
high value on
listening to and trying to understand each other,
we don't engage in interminable arguments and endless 
redoing of our
mathematics.  We are prepared to be convinced by others.   
Intellectually,
mathematics moves very quickly. Entire mathematical 
landscapes
change and change again in amazing ways during a single 
career.

When one considers how hard it is to write a computer 
program even approaching
the intellectual scope of a good mathematical paper, and 
how much greater
time and effort have to be put into it to make it 
``almost'' formally correct,
it is preposterous to claim that mathematics as we 
practice it is anywhere
near formally correct.

Mathematics as we practice it is much more formally 
complete and precise
than other sciences, but it is much less formally complete 
and precise
for its content than computer programs. 
The difference has to do not just with the amount
of effort: the kind of effort is qualitatively different.
In large computer programs, a tremendous proportion
of effort must be spent on myriad compatibility issues:
making sure that all definitions are consistent, 
developing ``good'' data
structures that have useful but not cumbersome generality,
deciding on the ``right'' generality for functions, {\it 
etc.}
The proportion of energy spent on the working part of a 
large program, as
distinguished from the bookkeeping part, is surprisingly 
small.
Because of compatibility issues that
almost inevitably escalate out of hand because the 
``right'' definitions
change as generality and functionality are added, computer 
programs
usually need to be rewritten frequently, often from scratch.

A very similar kind of effort would have to go into
mathematics to make it formally correct and complete.
It is not that formal correctness
is prohibitively difficult on a small scale---it's that
there are many possible choices of formalization
on small scales that translate to huge numbers of 
interdependent choices
in the large.  It is quite hard to make these choices 
compatible;
to do so would certainly entail going back and rewriting 
from scratch all
old mathematical papers whose results we depend on. 
It is also quite hard to come up with good technical choices
for formal definitions that will be valid in the variety 
of ways that
mathematicians want to use them and that will anticipate 
future
extensions of mathematics.  If we were to continue to 
cooperate,
much of our time would be spent with international standards
commissions to establish uniform definitions and resolve 
huge controversies.

Mathematicians can and do fill in gaps,
correct errors, and supply more detail and more careful 
scholarship when
they are called on or motivated to do so.  Our system is 
quite good at
producing reliable theorems that can be solidly backed up.
It's just that the reliability does not primarily
come from mathematicians formally checking formal 
arguments; it
comes from mathematicians thinking carefully and 
critically about
mathematical ideas.

\medskip
On the most fundamental level, the foundations of
mathematics are much shakier than the mathematics that we 
do.
Most mathematicians adhere to foundational principles that 
are
known to be polite fictions.  For example, it is a theorem 
that
there does not exist any way to ever actually construct or 
even define
a well-ordering of the real numbers.
There is considerable evidence (but no proof) that we can 
get away
with these polite fictions without being caught out, but 
that doesn't
make them right.  Set theorists construct many alternate
and mutually contradictory ``mathematical universes'' such 
that if
one is consistent, the others are too.  This leaves very 
little
confidence that one or the other is the right
choice or the natural choice.   G\"odel's incompleteness 
theorem
implies that there can be
no formal system that is consistent, yet powerful enough 
to serve as a
basis for all of the mathematics that we do.

In contrast to humans, computers are good at performing 
formal processes.
There are people working hard on the project of actually 
formalizing
parts of mathematics by computer, with actual formally 
correct formal
deductions.  I think this is a very big but very 
worthwhile project,
and I am confident that we will learn a lot from it.  The 
process will help
simplify and clarify mathematics.  In not too many years, 
I expect
that we will have interactive computer programs that can 
help people
compile significant chunks of formally complete and 
correct mathematics
(based on a few perhaps shaky but at least explicit 
assumptions),
and that they will become part of the standard 
mathematician's working
environment.

However, we should recognize that the humanly 
understandable and humanly
checkable proofs that we actually do are what is most 
important to us,
and that they are quite different from formal proofs.  For 
the present,
formal proofs are out of reach and mostly irrelevant:  we 
have good
human processes for checking mathematical validity.

\section{What motivates people to do mathematics?}

There is a real joy in doing mathematics, in learning ways 
of thinking
that explain and organize and simplify.  One can feel this 
joy discovering
new mathematics, rediscovering old mathematics, learning a 
way of
thinking from a person or text, or finding a new way to 
explain or to
view an old mathematical structure.

This inner motivation might lead us to think that we do 
mathematics
solely for its own sake.  That's not true:
the social setting is extremely important.  We are 
inspired by other
people, we seek appreciation by other people, and we like 
to help
other people solve their mathematical problems. What we 
enjoy changes
in response to other people.  Social interaction
occurs through face-to-face meetings.  It also occurs 
through
written and electronic correspondence, preprints, and 
journal articles.
One effect of this highly social system of mathematics is 
the tendency
of mathematicians to follow fads.  For the purpose of 
producing new
mathematical theorems this
is probably not very efficient: we'd seem to be better off 
having
mathematicians cover the intellectual field much more 
evenly.
But most mathematicians don't
like to be lonely, and they have trouble staying excited 
about
a subject, even if they are personally making progress, 
unless they have
colleagues who share their excitement.

In addition to our inner motivation and our informal 
social motivation for
doing mathematics, we are driven by considerations of 
economics and status.
Mathematicians, like other academics, do a lot of judging 
and being judged.
Starting with grades, and continuing through letters of 
recommendation,
hiring decisions, promotion decisions, referees reports, 
invitations to
speak, prizes, \dots  we are involved in many ratings, in 
a fiercely
competitive system.

\medskip
Jaffe and Quinn analyze the motivation to do mathematics 
in terms of a
common currency that many mathematicians believe in: 
credit for theorems.

I think that our strong communal emphasis on 
theorem-credits has a
negative effect on mathematical progress.
If what we are accomplishing is advancing human 
understanding of mathematics,
then we would be much better off recognizing and valuing a 
far
broader range of activity.   The people who see the way to 
proving
theorems are doing it in the context of a mathematical 
community;
they are not doing it on their own.   They depend on
understanding of mathematics that they glean
from other mathematicians.   Once a theorem has been 
proven, the
mathematical community depends on the social network to 
distribute the
ideas to people who might use them further---the print 
medium is far
too obscure and cumbersome.   

Even if one takes the narrow view that what we are 
producing is theorems, the
team is important.  Soccer can serve as a metaphor.
There might only be one or two goals during a soccer game,
made by one or two persons.
That does not mean that the efforts of all the others are 
wasted.
We do not judge players on a soccer team only by whether 
they personally
make a goal; we judge the team by its function as a team.

In mathematics, it often happens that a group of 
mathematicians
advances with a certain collection of ideas.  There
are theorems in the path of these advances that will 
almost inevitably be
proven by one person or another.  Sometimes the group of 
mathematicians
can even anticipate what these theorems are likely to be.  
It is much
harder to predict who will
actually prove the theorem, although there are usually a 
few ``point people''
who are more likely to score.  However, they are in a 
position to prove those
theorems because of the collective efforts of the team.  
The team has
a further function, in absorbing and making use of the 
theorems once
they are proven.   Even if one person could prove all the 
theorems in
the path single-handedly, they are wasted if nobody else 
learns them.

There is an interesting phenomenon concerning the 
``point'' people.
It regularly happens that someone who was in the middle of 
a pack
proves a theorem that receives wide recognition as being 
significant.
Their status in the community---their pecking 
order---rises immediately
and dramatically.  When this happens, they usually become
much more productive as a center of ideas and a source of 
theorems.
Why?   First, there is a large increase in self-esteem,
and an accompanying increase in productivity.
Second, when their status increases, people are more in 
the center of the
network of ideas---others take them more seriously.  
Finally and perhaps
most importantly, a mathematical breakthrough usually 
represents a new
way of thinking, and effective ways of thinking can 
usually be applied
in more than one situation.

This phenomenon convinces me that the entire mathematical 
community
would become much more productive if we open our eyes to 
the real values
in what we are doing.
Jaffe and Quinn propose a system of recognized roles 
divided into
``speculation'' and ``proving''.   Such a division only 
perpetuates
the myth that our progress is measured in units of 
standard theorems deduced.
This is a bit like the fallacy of the person who makes a 
printout of the
first 10,000 primes.  What we are producing is human 
understanding.
We have many different ways to understand and many 
different processes that
contribute to our understanding.  We will be more satisfied,
more productive and happier if we recognize
and focus on this.

\section{Some personal experiences}

Since this essay grew out of reflection on the misfit 
between my
experiences and the description of Jaffe and Quinn's,
I will discuss two personal experiences, including the one 
they alluded to.

I feel some awkwardness in this, because I do have regrets 
about aspects
of my career: if I were to do things over again with the
benefit of my present insights about myself
and about the process of mathematics, there is a lot that 
I would
hope to do differently. I hope that by describing these 
experiences
rather openly as I remember and understand them, I can 
help others
understand the process better and learn in advance.

\bigskip

First I will discuss briefly the theory of foliations, 
which was
my first subject, starting when I was a graduate student. 
(It doesn't
matter here whether you know what foliations are.)

At that time, foliations had become
a big center of attention among geometric topologists, 
dynamical systems people, and differential geometers.   I 
fairly rapidly
proved some dramatic theorems.   I proved a classification 
theorem for
foliations, giving a necessary and sufficient condition 
for a manifold to admit
a foliation.  I proved a number of other significant 
theorems.
I wrote respectable papers and published at least the
most important theorems.  It was hard to find the time to 
write to keep up
with what I could prove, and I built up a backlog.

An interesting phenomenon occurred.  Within a couple of 
years, a dramatic
evacuation of the field started to take place.  I heard 
from a number
of mathematicians that they were giving or receiving 
advice not to go into
foliations---they were saying that Thurston was cleaning 
it out.
People told me (not as a complaint, but as a compliment)
that I was killing the field.  Graduate students stopped 
studying
foliations, and fairly soon, I turned to other interests 
as well.

I do not think that the evacuation occurred because the 
territory was
intellectually exhausted---there were (and still are)
many interesting questions that remain and that are 
probably approachable.
Since those years, there have been interesting 
developments carried
out by the few people who stayed in the field or who 
entered the field,
and there have also been important developments in 
neighboring areas that
I think would have been much accelerated had 
mathematicians continued to
pursue foliation theory vigorously.

Today, I think there are few mathematicians who understand 
anything
approaching the state of the art of foliations as it lived 
at that time,
although there are some parts of the theory of foliations, 
including
developments since that time, that are still thriving.

I believe that two ecological effects were much more
important in putting a damper on the subject than any
exhaustion of intellectual resources that occurred.

First, the results I proved (as well as some important 
results of other people)
were documented in a conventional, formidable 
mathematician's style.  They depended heavily on readers 
who shared
certain background and certain insights.  The theory of 
foliations was a
young, opportunistic subfield, and the background was not 
standardized.
I did not hesitate to draw on any of the mathematics I had 
learned from others.
The papers I wrote did not (and could not) spend much time 
explaining
the background culture.   They documented top-level 
reasoning and
conclusions that I often had achieved after much 
reflection and effort.
I also threw out prize cryptic tidbits of insight,
such as ``the Godbillon-Vey invariant measures the helical 
wobble of a
foliation'', that remained mysterious to most 
mathematicans who read them. 
This created a high entry barrier:  I think many graduate 
students
and mathematicians were discouraged that it was hard to
learn and understand the proofs of key theorems.

Second is the issue of what is in it for other people in 
the subfield.
When I started working on foliations, I had the conception 
that what
people wanted was to know the answers.  I thought that 
what they sought was
a collection of powerful proven theorems that might be 
applied to answer
further mathematical questions.  But that's only one part 
of the story.
More than the knowledge, people want {\it personal 
understanding}.
And in our credit-driven system, they also want and need 
{\it theorem-credits}.

\bigskip

I'll skip ahead a few years, to the subject that Jaffe and 
Quinn alluded
to, when I began studying $3$-dimensional
manifolds and their relationship to hyperbolic geometry.  
(Again, it
matters little if you know what this is about.)
I gradually built up over a number of years a certain 
intuition for hyperbolic
three-manifolds, with a repertoire of constructions, 
examples and proofs.
(This process actually started when I was an 
undergraduate, and was strongly
bolstered by applications to foliations.) After a while,
I conjectured or speculated that all three-manifolds have 
a certain
geometric structure; this conjecture eventually became 
known as the
geometrization conjecture.  About two or three years later,
I proved the geometrization theorem
for Haken manifolds.  It was a hard theorem,
and I spent a tremendous amount of effort thinking
about it.  When I completed the proof, I spent a lot more 
effort checking
the proof, searching for difficulties and testing it 
against independent
information.

I'd like to spell out more what I mean when I say I proved 
this theorem.
It meant that I had a clear and complete flow of ideas, 
including details,
that withstood a great deal of scrutiny by myself and by 
others.
Mathematicians have many different styles of thought.
My style is not one of making broad sweeping but careless 
generalities,
which are merely hints or inspirations: I make clear 
mental models, and
I think things through.  My proofs have turned out to be 
quite reliable.
I have not had trouble backing up claims or producing 
details for things
I have proven.  I am good in detecting flaws in my own 
reasoning as well
as in the reasoning of others.

However, there is sometimes a huge expansion factor in 
translating
from the encoding in my own thinking to something that can 
be conveyed
to someone else.  My mathematical education was rather 
independent and
idiosyncratic, where for a number of
years I learned things on my own, developing personal 
mental models
for how to think about mathematics. This has often been a 
big advantage for
me in thinking about mathematics, because it's easy to 
pick up later the
standard mental models shared by groups of mathematicians.
This means that some concepts that I use freely and 
naturally in my personal
thinking are foreign to most mathematicians I talk to.  My 
personal
mental models and structures are similar in character to the
kinds of models groups of mathematicians share---but they 
are often
different models.  At the time of the formulation of the
geometrization conjecture, my understanding of hyperbolic 
geometry was a
good example.  A random continuing example is an 
understanding of
finite topological spaces, an oddball topic that can lend 
good
insight to a variety of questions but that is generally 
not worth
developing in any one case because there are standard
circumlocutions that avoid it.

Neither the geometrization conjecture nor its proof for 
Haken manifolds was
in the path of any group of mathematicians at the 
time---it went against
the trends in topology for the preceding 30 years, and it 
took
people by surprise.  To most topologists
at the time, hyperbolic geometry was an arcane side branch 
of
mathematics, although there were other groups of 
mathematicians
such as differential geometers who did understand it from 
certain
points of view. It took topologists a while just to 
understand what the
geometrization conjecture meant, what it was good for, and 
why it was relevant.

At the same time, I started writing notes
on the geometry and topology of $3$-manifolds, in 
conjunction with
the graduate course I was teaching.  I distributed them to 
a few people, and
before long many others from around the world were writing 
for copies.
The mailing list grew to about 1200 people to whom I was
sending notes every couple of months.    I tried to 
communicate my
real thoughts in these notes.  People ran many seminars
based on my notes, and I got lots of feedback. 
Overwhelmingly, the feedback ran something like ``Your 
notes are
really inspiring and beautiful, but I have to tell you 
that we spent
3 weeks in our seminar working out the details of \S\<$n.n$.
More explanation would sure help.''

I also gave many presentations to groups of mathematicians
about the ideas of studying $3$-manifolds
from the point of view of geometry, and about the proof of 
the geometrization
conjecture for Haken manifolds.  At the beginning, this 
subject was
foreign to almost everyone.  It was hard to 
communicate---the
infrastructure was in my head, not in the mathematical 
community.
There were several mathematical theories that fed into the 
cluster of
ideas: three-manifold topology, Kleinian groups, dynamical 
systems,
geometric topology, discrete subgroups of Lie groups, 
foliations,
Teichm\"uller spaces, pseudo-Anosov diffeomorphisms, 
geometric group theory,
as well as hyperbolic geometry.  

We held an AMS summer workshop at Bowdoin in 1980, where 
many mathematicans
in the subfields of low-dimensional topology, dynamical 
systems and
Kleinian groups came.

It was an interesting experience exchanging cultures.
It became dramatically clear how much proofs depend on the 
audience.
We prove things in a social context and address them to a 
certain audience. 
Parts of this proof I could communicate in two minutes
to the topologists, but the analysts would need an hour 
lecture
before they would begin to understand it.  Similarly, 
there were
some things that could be said in two minutes to the 
analysts that would
take an hour before the topologists would begin to get it. 
 And there
were many other parts of the proof which should take two 
minutes in
the abstract, but that none of the audience at the time 
had the
mental infrastructure to get in less than an hour.

At that time, there was practically no infrastructure
and practically no context for this theorem, so the 
expansion
from how an idea was keyed in my head to what I had to say 
to get it
across, not to mention how much energy the audience had to 
devote to
understand it, was very dramatic.

In reaction to my experience with foliations and in 
response to
social pressures, I concentrated most of my attention on 
developing and
presenting the infrastructure in what I wrote and in what 
I talked to
people about.  I explained the details to the few people who
were ``up'' for it.  I wrote some papers giving the 
substantive parts of the
proof of the geometrization theorem for Haken 
manifolds---for these papers,
I got almost no feedback.  Similarly,
few people actually worked through the harder and deeper
sections of my notes until much later.

The result has been that now quite a number of 
mathematicians have what
was dramatically lacking in the beginning:
a working understanding of the concepts and the 
infrastructure that are
natural for this subject.  There has been and there 
continues to be
a great deal of thriving mathematical activity.  By 
concentrating
on building the infrastructure and explaining and 
publishing definitions
and ways of thinking but being slow in stating or in 
publishing proofs
of all the ``theorems'' I knew how to prove, I left room
for many other people to pick up credit.  There has been 
room for
people to discover and publish other proofs of the 
geometrization theorem.
These proofs helped develop mathematical
concepts which are quite interesting in themselves, and 
lead to further
mathematics. 

What mathematicians most wanted and needed from me was to 
learn my ways of
thinking, and not in fact to learn my
proof of the geometrization conjecture for Haken manifolds.
It is unlikely that the proof of the general 
geometrization conjecture
will consist of pushing the same proof further.

\bigskip

A further issue is that people sometimes need or want an 
accepted and
validated result not in order to learn it, but so that they
can quote it and rely on it.

Mathematicians were actually very quick to accept my 
proof, and to start
quoting it and using it based on
what documentation there was, based on their experience 
and belief in me,
and based on acceptance by opinions of experts with whom I 
spent a lot of time
communicating the proof.  The theorem now is documented,
through published sources authored by me and by others,
so most people feel secure in quoting it;
people in the field certainly have not challenged me about 
its
validity, or expressed to me a need for details that are 
not available.

Not all proofs have an identical role in the logical 
scaffolding we are
building for mathematics. This particular proof probably 
has only
temporary logical value, although it has a high 
motivational value in
helping support a certain vision for the structure of 
$3$-manifolds.
The full geometrization conjecture is still a conjecture. 
It has been proven
for many cases, and is supported by a great deal of 
computer evidence
as well, but it has not been proven in generality.  I am 
convinced
that the general proof will be discovered; I hope before 
too many more
years.  At that point, proofs of special cases are likely 
to become obsolete.

Meanwhile, people who want to use the geometric technology 
are better off
to start off with the assumption
``Let $M^3$ be a manifold that admits a geometric 
decomposition,''
since this is more general than ``Let $M^3$ be a Haken 
manifold.''
People who don't want to use the technology or who are 
suspicious of it
can avoid it.  Even when a theorem about Haken manifolds 
can be proven using
geometric techniques, there is a high value in finding 
purely
topological techniques to prove it.

In this episode (which still continues) I think I have 
managed to avoid
the two worst possible outcomes: either for me not to let 
on that I
discovered what I discovered and proved what I proved, 
keeping it to
myself (perhaps with the hope of proving the Poincar\'e 
conjecture),
or for me to present an unassailable and hard-to-learn 
theory with no
practitioners to keep it alive and to make it grow.

\bigskip

I can easily name regrets about my career. I have not 
published
as much as I should. There are a number of mathematical 
projects
in addition to the geometrization theorem for Haken 
manifolds that
I have not delivered well or at all to the mathematical 
public.
When I concentrated more on developing the infrastructure 
rather than
the top-level theorems in the geometric theory of 
3-manifolds, I became
somewhat disengaged as the subject continued to evolve; 
and I have not
actively or effectively promoted the field or the careers 
of the excellent
people in it. (But some degree of disengagement seems to 
me an almost
inevitable by-product of the mentoring of graduate 
students and others:
in order to really turn genuine research directions over 
to others,
it's necessary to really let go and stop oneself from 
thinking about them
very hard.)

On the other hand, I have been busy and productive, in 
many different
activities. Our system does not create extra time for 
people like me to
spend on writing and research;
instead, it inundates us with many requests and 
opportunities for extra work,
and my gut reaction has been to say `yes' to many of these 
requests and
opportunities.  I have put a lot of effort into 
non-credit-producing activities
that I value just as I value proving theorems:   
mathematical
politics, revision of my notes into a book with a high 
standard of
communication, exploration of computing in mathematics,
mathematical education, development of new forms for 
communication of
mathematics through the Geometry Center
(such as our first experiment, the ``Not Knot'' video),
directing MSRI, {\it etc.}

\bigskip

I think that what I have done has not maximized my 
``credits''.
I have been in a position not to feel a strong need to 
compete for more
credits. Indeed, I began to feel strong challenges from 
other
things besides proving new theorems.

I do think that my actions have done well in stimulating 
mathematics.
\end{document}